\documentclass[12pt]{article}
\usepackage{amsmath,amssymb,amsthm}

\newtheorem*{thm}{Theorem}

\theoremstyle{definition}
\newtheorem*{abs}{ABSTRACT}
\newcommand{\ssection}[1]{%
  \section[#1]{\raggedright\sl #1}}

\begin{document}
\title{\textsc{A Simple Proof of Bohr's Inequality}}
\date{}
\author{\textsc{Vern I. Paulsen and Dinesh Singh}}
\maketitle

\begin{abs}
The classical inequality of Bohr concerning Taylor coeficients of bounded holomorphic functions on the unit disk, has proved to be of significance in answering in the negative the conjecture that if the non-unital von Neumann inequality held for a Banach algebra then it was necessarily an operator algebra. Here we provide a rather short and easy proof of the inequality. 
\end{abs}

\noindent
\ssection{Introduction}
\setcounter{footnote}{1}
\footnotetext{2000 \emph{Mathematics Subject Classification} Primary 46J10, 46J15; Secondary 30B10}
\setcounter{footnote}{2}
\footnotetext{\emph{Keywords and Phrases} Bohr's inequality}
Harald Bohr in 1914 obtained the following inequality in \cite{boh14}:
\begin{thm}
Let $f$ be in the disk algebra $A$ with Taylor expansion $f(z)=\sum_{n=0}^{\infty}a_nz^n$. Then $\sum_{n=0}^{\infty}|a_n|r^n\leq\|f\|_{\infty}$ for $0\le r\le \frac{1}{3}$ and $\frac{1}{3}$ is the best possible constant.
\end{thm}

Here $z=re^{i\theta}$, and $A$ is the Banach algebra of holomorphic functions on the open unit disk that have continuous boundary values and $\|f\|_{\infty}=sup\{|f(z)|:~|z|\le1\}$. Actually, Bohr obtained the inequality for $0\le r\le\frac{1}{6}$. Wiener, M. Riesz and Schur gave independent proofs and each established that $\frac{1}{3}$ is the best possible constant. See \cite{aiz00}, \cite{boa97}, \cite{boa00} and \cite{dix95}. Other proofs were also given by Sidon \cite{sid27} and by Tomic \cite{tom62}.

Operator algebraists began to take an interest in Bohr's inequality when Dixon \cite{dix95} used it to settle in the negative the conjecture that if the non-unital von Neumann inequality held for a Banach algebra then it was necessarily an operator algebra.

In \cite{pdspop}, we (along with Popescu) also established, using Bohr's inequality, that every Banach algebra has an equivalent norm that satisfies the non-unital von Neumann inequality thus providing large classes of Banach algebras that are not operator algebras and yet they satisfy the non-unital von Neumann inequality. In the same paper, Paulsen-Popescu-Singh also established versions of Bohr's inequality in the context of higher dimensional domains, for the non commutative disk algebra as well as for the reduced and full group $C^*-$algebras of the free group on $n$ generators. In \cite{din04} we provide a very general version of Bohr's inequality in the context of uniform algebras from which the classical inequality follows easily. In \cite{paul} we have extended Bohr's inequality to operator valued holomorphic functions as well as to the case of holomorphic functions on an annulus. For further connections see Aizenberg \cite{aiz00}, Boas and Khavinson \cite{boa97} and Boas \cite{boa00}.

\ssection{Preliminary observations}
In this note we provide a rather short and easy proof of Bohr's inequality using the well known property that the unit ball of the disk algebra $A$ is the closed convex hull of the set of finite Blaschke products, see [6, cor 2.4, page 196].

A finite Blaschke product is a function of the form $$\prod_{k=1}^{n}\dfrac{z-z_k}{1-\bar z_kz} $$ where $|z_k|<1$ for each $k$. We also need the following two easy and elementary observations: Setting $M(f,r)=\sum\limits_{n=0}^{\infty}|a_n|r^n$ where $f(z)=\sum\limits_{n=0}^{\infty}a_nz^n$, we see easily that
\begin{enumerate}
\item $M(\alpha f+\beta g,r)\leq |\alpha|M(f,r)+|\beta|M(g,r)$
\item $M(fg,r)\leq M(f,r)M(g,r)$\qquad$(0\leq r<1)$
\end{enumerate}
\ssection{Proof of the Theorem}
We assume with no loss of generality that $\|f\|_{\infty}=1$. Next, it is elementary to check that the inequality is true for a single Blaschke factor, $\varphi_a(z)=(z-a)/(1-\bar az)=-a+\sum_{n=1}^{\infty}\bar a^{n-1}(1-a\bar a)z^n.$ Indeed, $$M(\varphi_a, r) = \frac{|a|+(1-2|a|^2)r}{1-|a|r} \le 1,$$ for $0 \le r \le \frac{1}{3}$ and $|a|<1$, as a little calculus exercise shows. Next we note that if $f(z)=\sum_{n=0}^{\infty}a_nz^n$ and $g(z)=\sum_{n=0}^{\infty}b_nz^n$ are in $A$ then
$||a_n|-|b_n||\leq|a_n-b_n|=\left|\frac{1}{2\pi}\int_0^{2\pi}(f-g)e^{-in\theta}d\theta\right|\leq\|f-g\|_{\infty}$. Thus for $r\leq \frac{1}{3}$, $|M(f,r)-M(g,r)|\leq\frac{3}{2}\|f-g\|_{\infty}$. So in view of this inequality and by the norm density of the convex combinations of finite Blaschke factors in the unit ball of $A$ it is enough to establish Bohr's inequality for any convex combination of finite Blaschke products. However, by observation (1) above, it is enough to establish the inequality for a finite Blaschke product. But then again by observation (2), it suffices to establish Bohr's inequality for a single Blaschke product, which we have already done.  Hence the proof is complete. 

To see that $\frac{1}{3}$ is the best possible constant, verify that for $r> \frac{1}{3},$ the maximum of $M(\varphi_a, r)$ is attained for some $|a|<1$ and is strictly greater than $1.$ This is most easily seen by noting that the above formula is equal to $1$ at $|a|=1$ and is a decreasing function of $|a|$ at $1.$

The second author wishes to thank the Department of Mathematics, University of Houston and the Mathematical Sciences Foundation, Delhi, for their support.

\noindent
\textsc{Department of Pure Mathematics and Institute for Quantum Computing, University of Waterloo, Waterloo N2L 3G1}\\
\emph{E-mail address}: vpaulsen@uwaterloo.ca\\  
\textsc{Center for Lateral Innovation, Creativity and Knowledge, SGT University, Gurugram, Haryana 122001}\\
\emph{E-mail address}: dineshsingh1@gmail.com
\end{document}